\documentclass[11pt]{amsart}
\usepackage{latexsym}
\usepackage{fullpage}
\usepackage{amssymb}
\usepackage{amscd}
\usepackage{float}
\usepackage{graphicx}
\usepackage{amsfonts}
\usepackage{pb-diagram}
\usepackage{amsmath,amscd}
\usepackage{caption}
\usepackage{subcaption}
\usepackage{xcolor}

\newtheorem{thm}{Theorem}
\newtheorem{cor}{Corollary}

\newtheorem{defn}{Definition}
\newtheorem{remark}{Remark}

\begin{document}

\title[Pseudo links in handlebodies]
  {Pseudo links in handlebodies}

\author{Ioannis Diamantis}
\address{Department of Data Analytics and Digitalisation,
Maastricht University, School of Business and Economics,
P.O.Box 616, 6200 MD, Maastricht,
The Netherlands.
}
\email{i.diamantis@maastrichtuniversity.nl}


\setcounter{section}{-1}

\date{}

\begin{abstract}
In this paper we introduce and study the theory of pseudo links in the genus $g$ handlebody, $H_g$. Pseudo links are links with some missing crossing information that naturally generalize the notion of knot diagrams. The motivation for studying these relatively new knotted objects lies in the fact that pseudo links may be used to model DNA knots, since it is not uncommon for biologists to obtain DNA knots for which it is not possible to tell a positive from a negative crossing. We consider pseudo links in $H_g$ as mixed pseudo links in $S^3$ and we generalize the Kauffman bracket polynomial for the set of pseudo links in $H_g$. We then pass on the set of mixed pseudo braids, that is, pseudo braids whose closures are pseudo links in $H_g$, and we formulate the analogue of the Alexander theorem. It is worth mentioning that the theory of pseudo links is closely related to the theory of singular links and that all results in this paper may be used for studying singular links in $H_g$, by considering the pseudo knot theory as a quotient of the singular knot theory by one extra isotopy move.

\smallbreak
\bigbreak


\end{abstract}

\maketitle

\section{Introduction}\label{intro}

Pseudo diagrams of knots were introduce by Hanaki in \cite{H} as projections on the 2-sphere with over/under information at some of the double points. In particular, pseudo knots are standard knots whose projections contain some special crossings, called {\it pre-crossings}, that are considered neither positive, nor negative. The theory of {\it pseudo knots} is obtained by considering equivalence classes of pseudo diagrams under equivalence relations generated by a specific set of Reidemeister moves. With the use of the analogue of the Alexander theorem (\cite{BJW, D}), one may pass from pseudo links to pseudo braids, that is, classical braids with some crossing information missing. Moreover, in \cite{BJW}, the analogue of the Markov theorem for pseudo braid equivalence is presented and a sharpened version of this theorem is presented in \cite{D} with the use of the $L$-moves for pseudo braids. Finally, in \cite{HD}, the pseudo bracket polynomial, $<;>$, is defined for pseudo links in $S^3$, extending the Kauffman bracket polynomial for classical knots that is presented in \cite{LK}.

\smallbreak

In \cite{D5} the theory of pseudo links is generalized for the case of the Solid Torus, ST. More precisely, pseudo links in ST are presented as mixed links in $S^3$, and the pseudo bracket polynomial for pseudo links in $S^3$ is generalized for the case pseudo links in ST. Moreover, the mixed pseudo braid monoid of type B is introduced and studied in \cite{D5}, with the use of which, the analogues of the Alexander and the Markov theorems for pseudo links in ST are obtained. The theory of pseudo links in ST plays an important role for studying pseudo links in a genus $g$ handlebody, $H_g$, since $H_g$ is an oriented 3-manifold which is homeomorphic to a connected sum of $g$-solid tori.

\smallbreak

In this paper we extend the results of \cite{D5} for pseudo links in the genus $g$ handlebody, $H_g$. More precisely, we translate isotopy moves for pseudo links in $H_g$ to isotopy moves for mixed links in $S^3$ and we generalize the pseudo bracket polynomial for pseudo links in $H_g$, by first extending it to the genus two handlebody, $H_2$. We finally introduce and study the appropriate mixed pseudo braid monoid for this theory and we formulate the analogue of the Alexander theorem for pseudo links in $H_g$. It is worth mentioning that the knot theory in handlebodies plays an important role in the study of knot theories in 3-manifolds in general due to their connection with Heegaard splittings.

\bigbreak

The paper is organized as follows: In \S~\ref{prel} we recall results on pseudo links in $S^3$ and in the Solid Torus from \cite{BJW, HJMR, D, D5, La, LR1}. More precisely, in \S~\ref{secpl} we recall the definition of pseudo links in $S^3$ and we present the pseudo bracket polynomial for pseudo links in $S^3$. We also recall the definition of pseudo braids and the definition of the pseudo braid monoid and we state the analogues of the Alexander and the Markov theorems. Moreover, we state a sharpened version of the analogue of the Markov theorem for pseudo links in $S^3$ with the use of the $L$-moves defined in this topological setting. In \S~\ref{ktst} we present similar results for the category of pseudo links in the solid torus, ST. We first view pseudo links in ST as mixed pseudo links in $S^3$ and we present the analogue of the Reidmeister theorem for pseudo links in ST. We then present the definition of the pseudo bracket polynomial for pseudo links in ST and we then pass on the category of mixed pseudo braids, that is, pseudo braids in ST whose closures are pseudo links in ST. We state the analogue of the Alexander theorem for pseudo links in ST and with the use of the mixed pseudo braid monoid, we present the analogue of the Markov theorem for pseudo links in ST. In \S~\ref{tpl} we introduce and study the theory of pseudo links and pseudo braids in the handlebody of genus $g$, $H_g$. In particular, in \S~\ref{plhg1} we set up the topological setting in order to extend the definition of the pseudo bracket polynomial for pseudo links in $H_g$ (see \S~\ref{plhg2}). In \S~\ref{plhg3}, we introduce the geometric mixed pseudo braids in $S^3$, that represent mixed links in $H_g$, and with the use of the $L$-moves, we present a braiding algorithm for pseudo links in $H_g$.

\section{Preliminaries}\label{prel}

\subsection{Pseudo links in $S^3$}\label{secpl}

In this subsection we present the theory of pseudo links and pseudo braids in $S^3$. In particular, we first present equivalence moves for pseudo links in $S^3$ and we recall the definition of the pseudo bracket polynomial that generalizes the Kauffman bracket polynomial for classical knots in $S^3$. We then pass on the level of pseudo braids and we present the analogues of the Alexander and the Markov theorems for pseudo links in $S^3$, with the use of the pseudo braid monoid $PM_n$.

\smallbreak

\begin{figure}[ht]
\begin{center}
\includegraphics[width=1.6in]{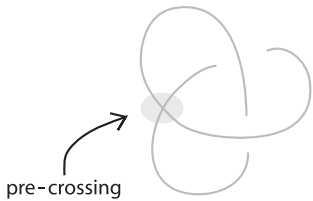}
\end{center}
\caption{A pseudo knot.}
\label{pk1}
\end{figure}

A {\it pseudo diagram} of a knot consists of a regular knot diagram with some missing crossing information, that is, there is no information about which strand passes over and which strand passes under the other. These undetermined crossings are called {\it pre-crossings} or {\it pseudo-crossings} (for an illustration see Figure~\ref{pk1}).

\begin{figure}[ht]
\begin{center}
\includegraphics[width=5.9in]{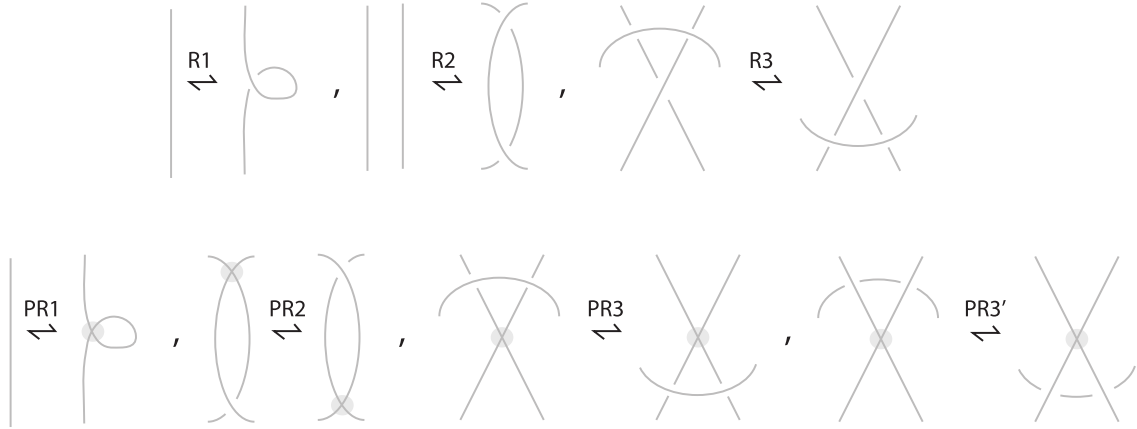}
\end{center}
\caption{Reidemeister moves for pseudo knots.}
\label{reid}
\end{figure}

\begin{defn}\rm
{\it Pseudo knots} are defined as equivalence classes of pseudo diagrams under an appropriate choice of Reidemeister moves that are illustrated in Figure~\ref{reid}. Moreover, a {\it pseudo link} is a collection of knotted pseudo knots.
\end{defn}

In \cite{HD} the pseudo bracket polynomial, $<\, ;>$, is defined for pseudo links in $S^3$, extending the Kauffman bracket polynomial for classical knots presented in \cite{LK}. The main difference between the pseudo and the classical bracket polynomial lies in the fact that in order to define skein relations on pre-crossings for the case of pseudo links, the orientation of a diagram is needed. More precisely, we have the following:

\begin{defn}\label{pkaufb}\rm
Let $L$ be an oriented pseudo link in $S^3$. The {\it pseudo bracket polynomial} of $L$ is defined by means of the following relations:
\begin{figure}[ht]
\begin{center}
\includegraphics[width=3.5in]{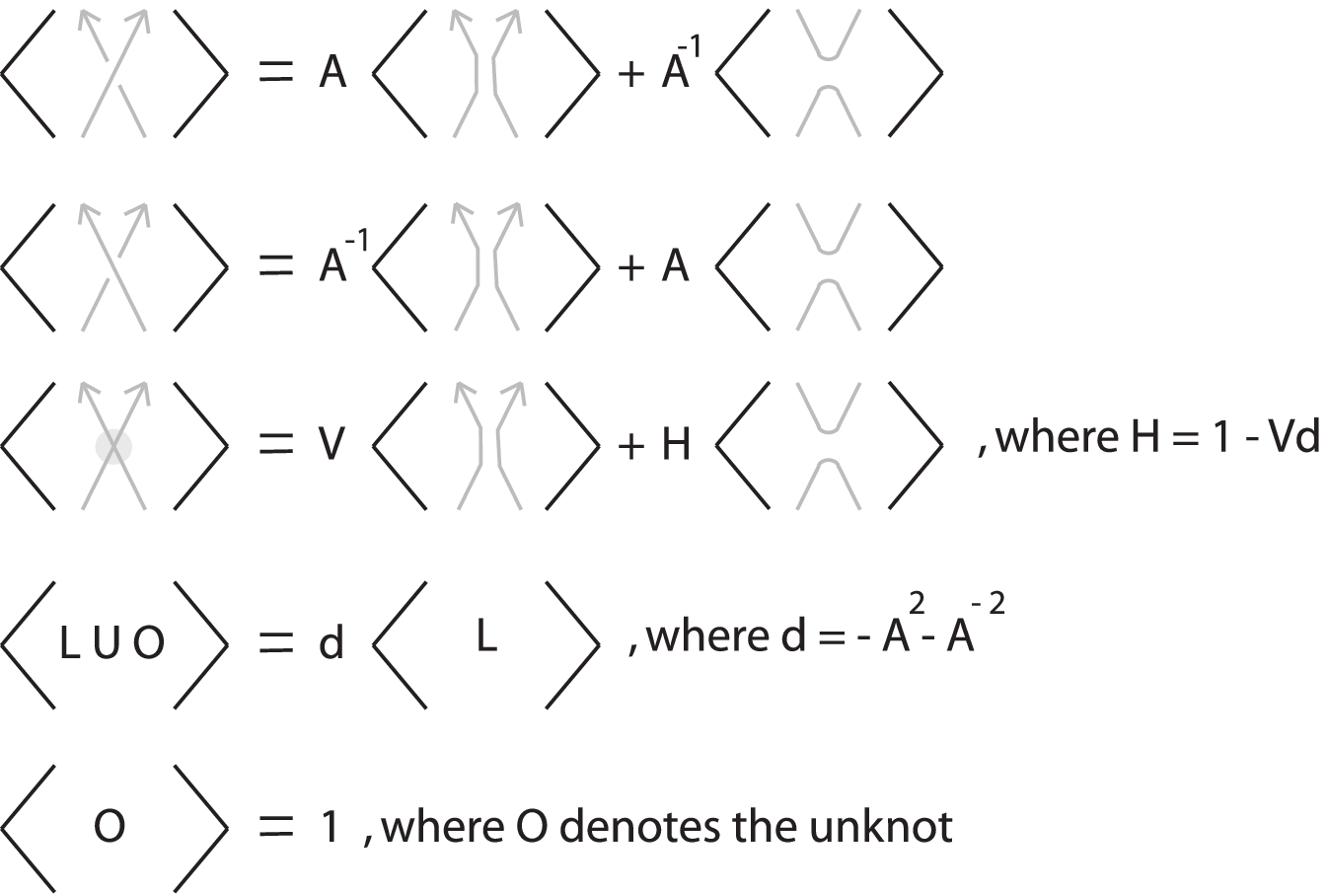}
\end{center}
\label{pkb}
\end{figure}
\end{defn}

The pseudo bracket is invariant under the Reidemeister moves 2 and 3 and under the pseudo moves PR1, PR2 and PR3 (\cite{HD} Theorem~1). As shown in \cite{LK}, in order to obtain an invariant for pseudo knots in $S^3$, the pseudo bracket polynomial is normalized using the writhe, leading to the {\it normalized pseudo bracket polynomial} (see Corollary~2 \cite{HD}).

\begin{thm}
Let $K$ be a pseudo diagram of a pseudo knot. The polynomial
\[
P_K(A, V)\ =\ (-A^{-3})^{w(K)}\, <K>,
\]
\noindent where $w(K):=\underset{c\in C(K)}{\sum}\, sgn(c)$, $C(k)$ the set of classical crossings of $K$ and $<K>$ the pseudo bracket polynomial of $K$, is an invariant of pseudo knots in $S^3$.
\end{thm}

We now introduce the pseudo braid monoid, $PM_n$, following \cite{BJW}.

\begin{figure}[ht]
\begin{center}
\includegraphics[width=2.6in]{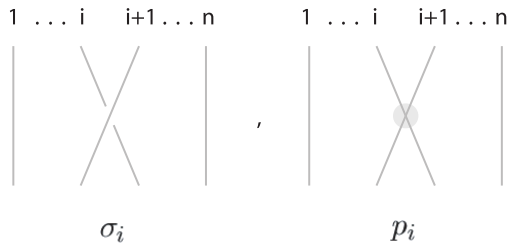}
\end{center}
\caption{The pseudo braid monoid generators.}
\label{gens}
\end{figure}

\begin{defn}\label{pmn}\rm
The monoid of pseudo braids, $PM_n$, is the monoid generated by $\sigma_i^{\pm 1}, p_i, i=1, \ldots, n-1$, illustrated in Figure~\ref{gens}, where $\sigma_i^{\pm 1}$ generate the braid group $B_n$ and $p_i$ satisfy the following relations:
\[
\begin{array}{rlcll}
i. & p_i\, p_j & = & p_j\, p_i, & {\rm if}\ |i-j|\geq 2\\
&&&&\\
ii. & p_i\, \sigma_j^{\pm 1} & = & \sigma_j^{\pm 1}\, p_i, & {\rm if}\ |i-j|\geq 2\\
&&&&\\
iii. & p_i\, \sigma_i^{\pm 1} & = & \sigma_i^{\pm 1}\, p_i, & i=1, \ldots, n-1\\
&&&&\\
iv. & \sigma_i\, \sigma_{i+1}\, p_i & = & p_{i+1}\, \sigma_i\, \sigma_{i+1}, & i=1, \ldots, n-2\\
&&&&\\
v. & \sigma_{i+1}\, \sigma_i\, p_{i+1} & = & p_{i}\, \sigma_{i+1}\, \sigma_i, & i=1, \ldots, n-2\\
\end{array}
\]
\end{defn}

\begin{remark}\label{propb2}\rm
It is worth mentioning that the monoid of pseudo braids is isomorphic to the {\it singular braid monoid}, $SM_n$, namely, the algebraic counterpart analogue of {\it singular knots}, i.e. knots that contain a finite number of self-intersections (Proposition 2.3 \cite{BJW}). Moreover, in \cite{FKR} it is shown that $SM_n$ embeds in a group, the singular braid group $SB_n$. It follows that $PM_n$ embeds in a group also, the pseudo braid group $PB_n$, generated by $\sigma_i$ and $p_i$, $i=1, \ldots, n-1$, satisfying the same relations as $PM_n$.
\end{remark}

Define now the {\it closure} of a pseudo braid as in the standard case. By considering $PM_n \subset PM_{n+1}$, we can consider the inductive limit $PM_{\infty}$. Using the analogue of the Alexander theorem for singular knots (\cite{Bi}), in \cite{BJW} the analogue of the Alexander theorem for pseudo links is presented. In particular:

\begin{thm}[{\bf Alexander's theorem for pseudo links}] \label{alexpl}
Every pseudo link can be obtained by closing a pseudo braid.
\end{thm}

We now state the analogue of the Markov theorem for pseudo links in $S^3$ (\cite{BJW}).

\begin{thm}[{\bf The analogue of the Markov Theorem for pseudo braids}] \label{markpl}
Two pseudo braids have isotopic closures if and only if one can be obtained from the other by a finite sequence of the following moves:
\[
\begin{array}{lllcll}
{Conjugation:} &  \alpha & \sim & \beta^{\pm 1}\, \alpha\, \beta^{\mp 1}, & {\rm for}\ \alpha \in PM_n\ \&\ \beta \in B_n,\\
&&&&\\
{Commuting:} &  \alpha\, \beta & \sim & \beta\, \alpha, & {\rm for}\ \alpha,\, \beta \in PM_n,\\
&&&&\\
{Stabilization:} &  \alpha & \sim & \alpha\, \sigma_n^{\pm 1}, & {\rm for}\ \alpha \in PM_n,\\
&&&&\\
{Pseudo-stabilization:} &  \alpha & \sim & \alpha\, p_n, & \alpha \in PM_n.\\
\end{array}
\]
\end{thm}

In \cite{D} a sharpened version of the analogue of the Markov theorem for pseudo braids is presented, with the use of the $L$-moves. $L$-moves make up an important tool for braid equivalence in any topological setting and they allow us to formulate sharpened versions of the analogues of the Markov theorems. For more details the reader is referred to \cite{La} and references therein.

\begin{figure}[ht]
\begin{center}
\includegraphics[width=4.8in]{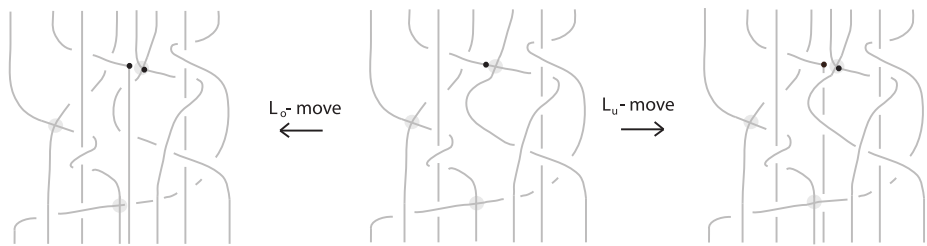}
\end{center}
\caption{$L$-moves for pseudo braids.}
\label{lm}
\end{figure}

\begin{defn}\label{lmdefn}\rm
An {\it $L$-move} on a braid $\beta$, consists in cutting an arc of $\beta$ open and pulling the upper cut-point downward and the lower upward, so as to create a new pair of braid strands with corresponding endpoints (on the vertical line of the cut-point), and such that both strands cross entirely {\it over} or {\it under} with the rest of the braid. Stretching the new strands over will give rise to an {\it $L_o$-move\/} and under to an {\it  $L_u$-move\/} as shown in Figure~\ref{lm} by ignoring all pre-crossings.
\end{defn}

Similarly, $L$-moves for pseudo braids are defined as follows: the two strands that appear after the performance of an $L$-move should cross the rest of the braid only with real crossings (all over in the case of an $L_o$-move or all under in the case of an $L_u$-move). For an illustration see Figure~\ref{lm}. 

We now recall the sharpened version of the analogue of the Markov theorem for pseudo braids using the $L$-moves (see Theorem~8 \cite{D}).

\begin{thm}[{\bf $L$-move braid equivalence for pseudo braids}] \label{markpll}
Two pseudo braids have isotopic closures if and only if one can be obtained from the other by a finite sequence of the following moves:
\[
\begin{array}{llcll}
{L{\text -}moves} &   &  &   & \\
&&&&\\
{Commuting:} &  \alpha\, \beta & \sim & \beta\, \alpha, & {\rm for}\ \alpha,\, \beta \in PM_n,\\
&&&&\\
{Pseudo{\text -}Stabilization:} &  \alpha & \sim & \alpha\, p_n, & \in PM_{n+1}.\\
\end{array}
\]
\end{thm}

\begin{remark}\label{Lrem}\rm
In \cite{D}, it is shown that $L$-moves on pseudo knots may be used in order to obtain the analogue of the Alexander theorem for pseudo knots following the braiding algorithm in \cite{LR1}.
\end{remark}

\subsection{Pseudo links in ST}\label{ktst}

We now view ST as the complement of a solid torus in $S^3$ and we present results from \cite{D5}. An oriented link $L$ in ST can be viewed as an oriented \textit{mixed link} in $S^{3}$, that is, a link in $S^{3}$ consisting of the unknotted fixed part $\widehat{I}$ that represents the complementary solid torus in $S^3$, and the moving part $L$ that links with $\widehat{I}$. Similarly, pseudo links in ST can be viewed as {\it mixed pseudo links} in $S^3$, that is, mixed links with some crossing information missing. Note that all crossings between the fixed and the moving parts of a mixed pseudo link are standard crossings and not pre-crossings.

\smallbreak

\begin{figure}[ht]
\begin{center}
\includegraphics[width=1.3in]{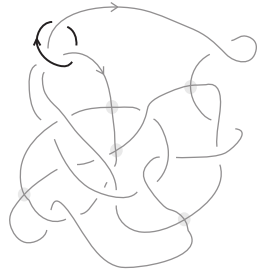}
\end{center}
\caption{A mixed pseudo link in $S^3$.}
\label{mplink}
\end{figure}

A \textit{mixed pseudo link diagram }is a pseudo diagram $\widehat{I}\cup \widetilde{L}$ of $\widehat{I}\cup L$ on the plane of $\widehat{I}$, where this plane is equipped with the top-to-bottom direction of $I$. For an illustration see Figure~\ref{mplink}.

\begin{figure}[ht]
\begin{center}
\includegraphics[width=4.9in]{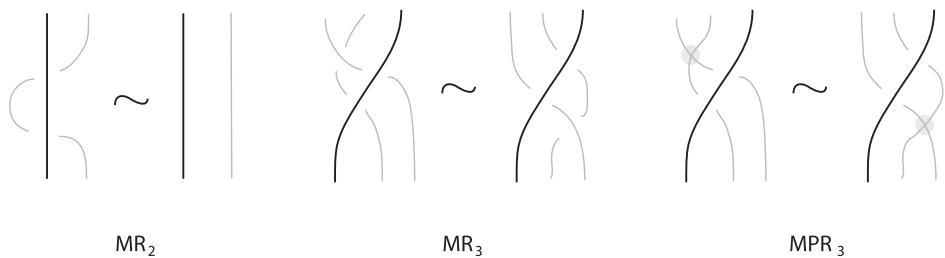}
\end{center}
\caption{Mixed Reidemeister moves.}
\label{mpr}
\end{figure}

In terms of pseudo diagrams, isotopy for pseudo links in ST is translated on the level of mixed pseudo links in $S^3$ by means of the following theorem:

\begin{thm}[{\bf Reidemeister's theorem for mixed pseudo links}]\label{reidplink}
Two mixed pseudo links in $S^3$ are isotopic if and only if they differ by a finite sequence of the classical and the pseudo Reidemeister moves illustrated in Figure~\ref{reid} for the standard part of the mixed pseudo links, and moves that involve the fixed and the standard part of the mixed pseudo links, that are called mixed Reidemeister moves, and which are illustrated in Figure~\ref{mpr}.
\end{thm}

Note that there are two additional mixed Reidemeister moves where the fixed part of the mixed braid lies below the moving part.

\smallbreak

\begin{figure}[H]
\begin{center}
\includegraphics[width=3.8in]{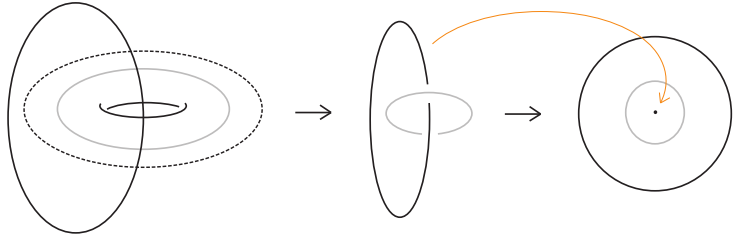}
\end{center}
\caption{ST as a punctured disk.}
\label{ppsttt}
\end{figure}

The pseudo bracket polynomial for pseudo links in $S^3$ can be naturally generalized for pseudo links in ST. For that reason, we view ST as a punctured disk, that is, a disk with a hole in its center, representing the complementary solid torus in $S^3$ (for an illustration see Figure~\ref{ppsttt}). 

\begin{defn}\label{pkaufbst}\rm
Let $L$ be an oriented pseudo link in ST. The {\it pseudo bracket polynomial} of $L$ is defined by means of the relations in Definition~\ref{pkaufb} together with the following relations:
\begin{figure}[H]
\begin{center}
\includegraphics[width=1.5in]{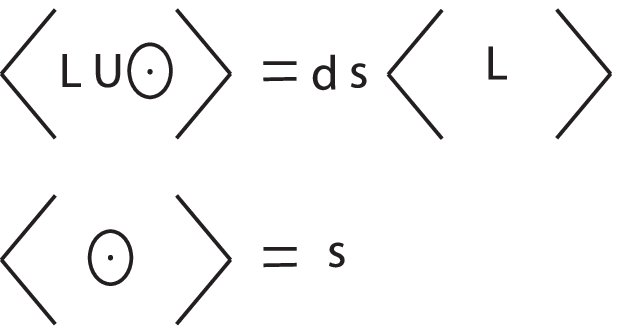}
\end{center}
\label{pkbst}
\end{figure}
\end{defn}

Obviously, the pseudo bracket polynomial for pseudo knots in ST is invariant under Reidemeister moves 2 and 3 and all pseudo moves, and similarly to the case of pseudo links in $S^3$, we may normalize the pseudo bracket polynomial for pseudo links in ST in order to obtain an invariant for pseudo links in ST. Namely:

\begin{cor}
Let $K$ be a pseudo diagram of a pseudo link in ST. The polynomial
\[
P_K(A, V, s)\ =\ (-A^{-3})^{w(K)}\, <K>,
\]
\noindent where $<K>$ denotes the pseudo bracket polynomial of $K$, is an invariant of pseudo links in ST.
\end{cor}

\begin{figure}[ht]
\begin{center}
\includegraphics[width=3.3in]{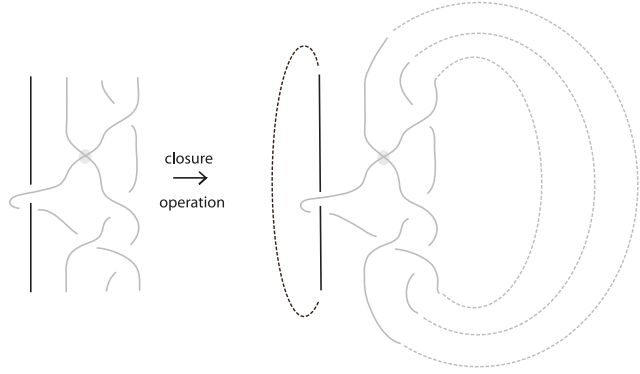}
\end{center}
\caption{A mixed pseudo braid and its closure to a mixed pseudo link.}
\label{cmpl}
\end{figure}

By the analogue of the Alexander theorem for pseudo links in the Solid Torus ST (\cite{D5}), a mixed pseudo link diagram $\widehat{I}\cup \widetilde{L}$ of $\widehat{I}\cup L$ may be turned into a \textit{mixed pseudo braid} $I\cup \beta $ with isotopic closure, where the {\it closure} of a mixed pseudo braid is defined as in the case of classical braids in $S^3$.  Mixed pseudo braids are pseudo braids in $S^{3}$ where, without loss of generality, their first strand represents $\widehat{I}$, the fixed part, and the other strands, $\beta$, represent the moving part $L$. The subbraid $\beta$ shall be called the \textit{moving part} of $I\cup \beta$. For an illustration see Figure~\ref{cmpl}.

Algebraic structures related to mixed pseudo braids are presented in \cite{D5}, by decomposing every isotopy in a sequence of elementary isotopies which correspond to the relations in Definition~\ref{mpbm}. As explained in \cite{D5}, these relations form a complete set of relations for the {\it mixed pseudo braid monoid} $PM_{1, n}$, the counterpart of the Artin's braid group of type B for pseudo links.

\begin{figure}[ht]
\begin{center}
\includegraphics[width=1in]{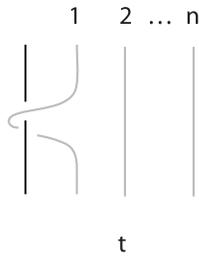}
\end{center}
\caption{The $t$-generator of $PM_{1,n}$.}
\label{genn}
\end{figure}

\begin{defn}\label{mpbm}\rm
The {\it mixed pseudo braid monoid of type B} $PM_{1, n}$ is defined as the monoid generated by the standard braid generators $\sigma_i^{\pm 1}$'s of $B_n$, the pseudo braid generators $p_{i}$'s of $PM_n$ and the looping generator $t$, illustrated in Figure~\ref{genn}, satisfying the following relations:

\[
\begin{array}{rcll}
\sigma_i\, \sigma_j & = & \sigma_j\, \sigma_i, & {\rm for}\ |i-j|>1,\\
&&&\\
\sigma_i\, \sigma_j\, \sigma_i & = & \sigma_j\, \sigma_i\, \sigma_j, & {\rm for}\ |i-j|=1,\\
&&&\\
p_i\, p_j & = & p_j\, p_i, & {\rm for}\ |i-j|>1,\\
&&&\\
\sigma_i\, p_j & = & p_j\, \sigma_i, & {\rm for}\ |i-j|>1,\\
&&&\\
\sigma_i\, \sigma_j\, p_i & = & p_j\, \sigma_i\, \sigma_j, & {\rm for}\ |i-j|=1,\\
&&&\\
t\, \sigma_i & = & \sigma_i\, t, & {\rm for}\ i>1,\\
&&&\\
t\, p_i & = & p_i\, t, & {\rm for}\ i>1,\\
&&&\\
t\, \sigma_1\, t\, \sigma_1 & = & \sigma_1\, t\, \sigma_1\, t, &\\
&&&\\
t\, \sigma_1\, t\, p_1 & = & p_1\, t\, \sigma_1\, t, &\\
&&&\\
\sigma_1\, \sigma_1^{-1} & = & t\, t^{-1}\ =\ 1.&
\end{array}
\]
\end{defn}

\begin{remark}\rm
It is worth mentioning that $PM_{1, n}$ is isomorphic to the {\it singular braid monoid of type B}, presented in \cite{V} (Theorem~7 \cite{D5}).
\end{remark}

We finally have that isotopy in ST is translated on the level of mixed pseudo braids by means of the following theorem:

\begin{thm}[{\bf The analogue of the Markov theorem for mixed pseudo braids}] \label{marktpb}
Two mixed pseudo braids have equivalent closures if and only if one can  obtained from the other by a finite sequence of the following moves:
\[
\begin{array}{llcll}
{\rm Commutation:} &  \alpha\, p_i & \sim & p_i\, \alpha, & {\rm for\ all}\ \alpha \in PM_{1, n},\\
&&&&\\
{\rm Conjugation:} & \beta & \sim & \alpha^{\pm 1}\, \beta\, \alpha^{\mp 1} & {\rm for\ all}\ \beta,\, \alpha\in PM_{1, n}\ {\rm such\ that}\\
&&&& \alpha\ {\rm contains\ only\ standard\ crossings},\\
&&&&\\
{\rm Real-Stabilization:} &  \alpha & \sim & \alpha\, \sigma_n^{\pm 1}, & {\rm for\ all}\ \alpha \in PM_{1, n},\\
&&&&\\
{\rm Pseudo-Stabilization:} &  \alpha & \sim & \alpha\, p_n, & {\rm for\ all}\ \alpha \in PM_{1, n}.\\
\end{array}
\]
\end{thm}

\section{Pseudo Links and pseudo braids in $H_g$}\label{tpl}

In this section we introduce and study the theory of pseudo links in the genus $g$ handlebody $H_g$. Throughout this section, we will be using results from \cite{OL}, where the classical knot theory of $H_g$ is introduced and studied.

\subsection{Pseudo links in $H_g$}\label{plhg1}

\begin{figure}[ht]
\begin{center}
\includegraphics[width=1.9in]{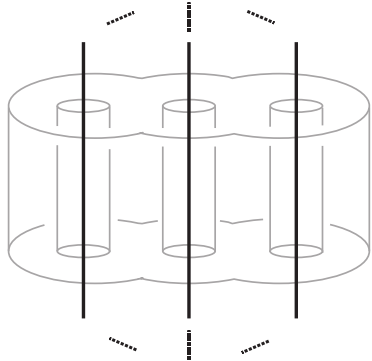}
\end{center}
\caption{The genus 3 handlebody.}
\label{h3}
\end{figure}

We consider $H_g$ to be $S^3\backslash \{{\rm open\ tubular\ neighborhood\ of}\ I_g \}$, where $I_g$ denotes the point-wise fixed identity braid on $g$ indefinitely extended strands meeting at the point at infinity. Thus $H_g$ may be represented in $S^3$ by the braid $I_g$ (for an illustration see Figure~\ref{h3}). 

\begin{figure}[ht]
  \begin{subfigure}[b]{0.28\textwidth}
    \includegraphics[width=\textwidth]{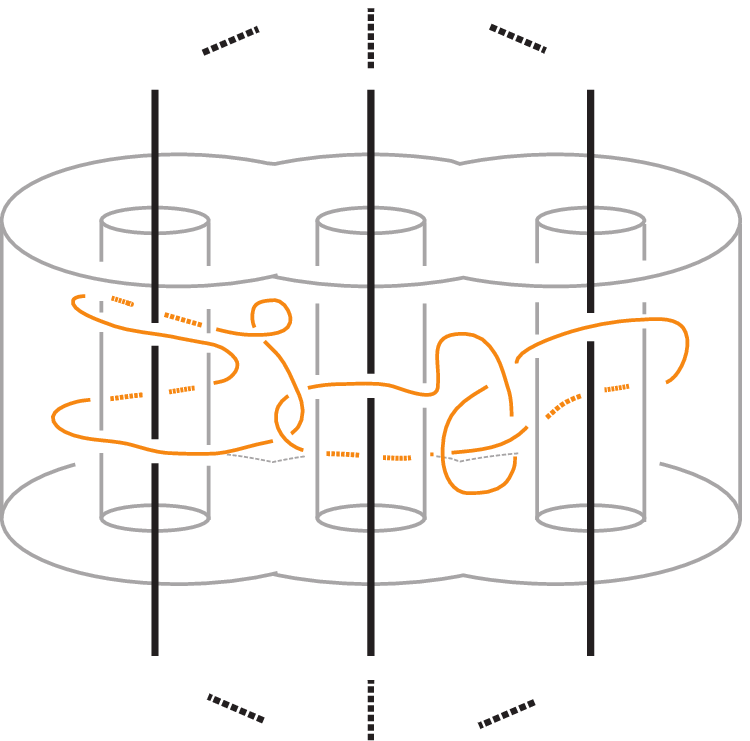}
    \caption{A mixed link in $H_3$.}
    \label{kHg}
  \end{subfigure}
  \begin{subfigure}[b]{0.24\textwidth}
    \includegraphics[width=\textwidth]{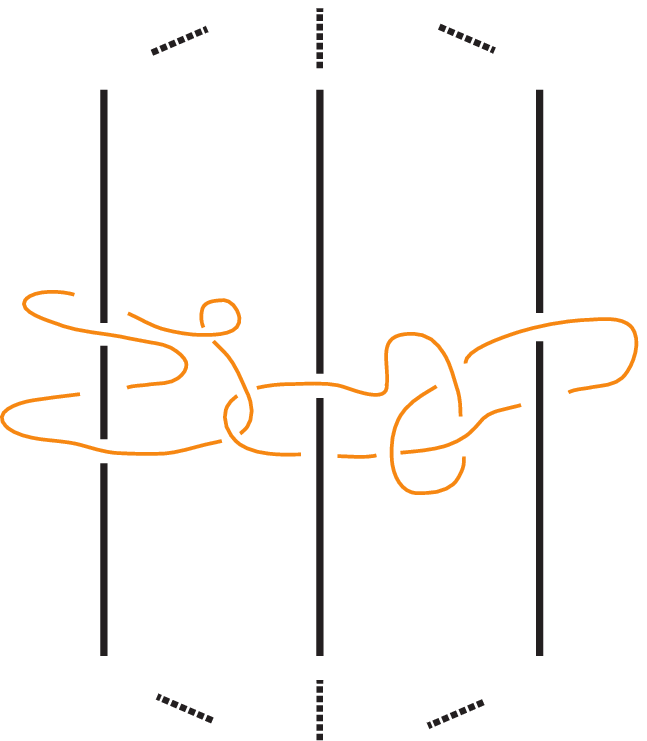}
    \caption{Representation of a mixed link in $H_3$.}
    \label{kHg2}
  \end{subfigure}
\end{figure}

An oriented link $L$ in $H_g$ can be represented by an oriented \textit{mixed link} in $S^{3}$, as illustrated in Figures~\ref{kHg} \& \ref{kHg2} (for more information the reader is referred to \cite{LR1, OL, DL1, D5}. A \textit{mixed link diagram} is a diagram $I_g\cup \widetilde{L}$ of $I_g\cup L$ on the plane of $I_g$, where this plane is equipped with the top-to-bottom direction of $I_g$. Similarly, an (oriented) pseudo link in $H_g$ can be represented as a {\it mixed pseudo link} in $S^3$, namely, a pseudo link in $S^{3}$ consisting of the fixed part $I_g$ and the moving part $L$ that links with $I_g$. Moreover, a \textit{mixed pseudo link diagram} is a diagram $I_g\cup \widetilde{L}$ of $I_g\cup L$ on the plane of $I_g$, where this plane is equipped with the top-to-bottom direction of $I_g$.

Isotopy for pseudo links in $H_g$ is then naturally translated in terms of isotopy for mixed pseudo links in $S^3$. In particular, two (oriented) pseudo links in $H_g$ are isotopic if and only if any two corresponding mixed pseudo link diagrams of theirs in $S^3$, differ by a finite sequence of the three Reidemeister moves, that involve classical crossings on their moving parts and the Mixed-Reidemeister moves, illustrated in Figure~\ref{mpr}.

\subsection{The pseudo bracket polynomial}\label{plhg2}

\begin{figure}[!ht]
\begin{center}
\includegraphics[width=1.9in]{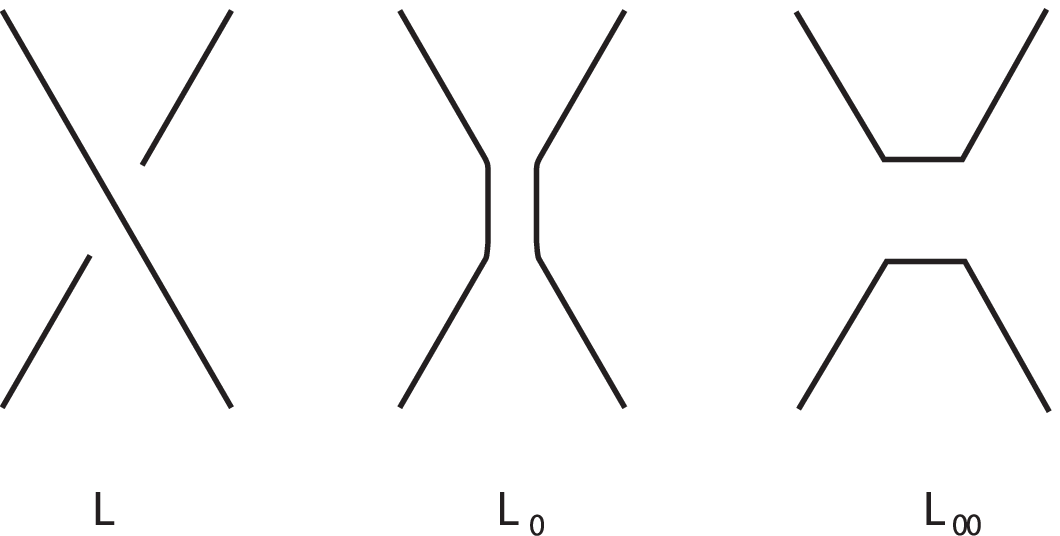}
\end{center}
\caption{The links $L$, $L_0$ and $L_{\infty}$ locally.}
\label{skein}
\end{figure}

In this subsection we extend the pseudo bracket polynomial for pseudo links in $H_g$. Our motivation comes from the theory of {\it skein modules}, which were independently introduced by Przytycki \cite{P} and Turaev \cite{Tu} as generalizations of knot polynomials in $S^3$ to knot polynomials in arbitrary 3-manifolds. In particular, let $M$ be an oriented $3$-manifold and $\mathcal{L}_{{\rm fr}}$ be the set of isotopy classes of unoriented framed links in $M$. Let $R=\mathbb{Z}[A^{\pm1}]$ be the Laurent polynomials in $A$ and let $R\mathcal{L}_{{\rm fr}}$ be the free $R$-module generated by $\mathcal{L}_{{\rm fr}}$. Let $\mathcal{S}$ be the ideal generated by the skein expressions $L-AL_{\infty}-A^{-1}L_{0}$ and $L \bigsqcup {\rm O} - (-A^2-A^{-1})L$, where $L_{\infty}$ and $L_{0}$ are exactly the same except inside a small ball, where the diagrams are illustrated in Figure~\ref{skein}. Note that blackboard framing is assumed. 

Then the {\it Kauffman bracket skein module} of $M$, KBSM$(M)$, is defined to be:

\begin{equation*}
{\rm KBSM} \left(M\right)={\raise0.7ex\hbox{$
R\mathcal{L}_{{\rm fr}} $}\!\mathord{\left/ {\vphantom {R\mathcal{L_{{\rm fr}}} {\mathcal{S} }}} \right. \kern-\nulldelimiterspace}\!\lower0.7ex\hbox{$ S  $}}.
\end{equation*}

\begin{figure}[!ht]
\begin{center}
\includegraphics[width=2.6in]{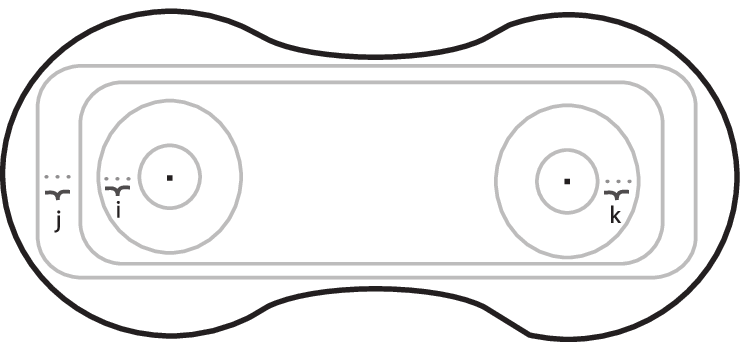}
\end{center}
\caption{The basis of KBSM($H_2$), $B_{H_2}$.}
\label{przt}
\end{figure}

In order to construct an invariant for pseudo links in $H_g$, we start by considering pseudo links in the handlebody of genus 2, $H_2$. In \cite{P} the Kauffman bracket skein module of the handlebody of genus 2, $H_2$, is computed using diagrammatic methods by means of the following theorem:

\begin{thm}[\cite{P}]\label{tprz}
The Kauffman bracket skein module of $H_2$, KBSM($H_2$), is freely generated by an infinite set of generators $\left\{x^i\, y^j\, z^k,\ (i, j, k)\in \mathbb{N}\times \mathbb{N}\times \mathbb{N}\right\}$, where $x^i\, y^j\, z^k$ is shown in Figure~\ref{przt}.
\end{thm}

\noindent Note that in \cite{P}, $H_2$ is represented as a twice punctured annulus and diagrams in $H_2$ determine a framing of links in $H_2$. The two dots in Figure~\ref{przt} represent the punctures. For more details the reader is referred to \cite{B}.

\begin{figure}[!ht]
\begin{center}
\includegraphics[width=1.2in]{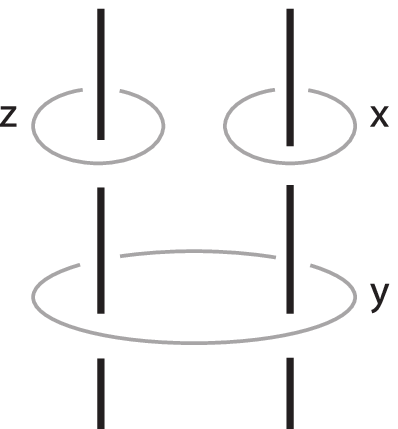}
\end{center}
\caption{The basic elements of KBSM($H_2$) as mixed links.}
\label{przt2}
\end{figure}

In our setting, the element $x^k\, z^l\, y^m$ denotes $k$-copies of the closed curve $x$, $l$-copies of the closed curve $z$ and $m$-copies of the closed curve $y$, as illustrated in Figure~\ref{przt2}. Note also that $x^0\, z^0\, y^0$ represents the unknot.

Equivalently, the above states that applying the Kauffman bracket skein relations on links in $H_2$, we end up with a finite formal sum where the terms consist of elements of the form $x^k\, z^l\, y^m$, that is, a number of unknots, and a number of $x, y$ and $z$ curves that we need to ``normalize'' in order to obtain a polynomial invariant for links in $H_2$. Normalizing these curves means to assign a specific variable to each one on the level of KBSM($H_2$). Similarly, we would have to normalize additional curves for obtaining the analogue of the Kauffman bracket polynomial for links in $H_g$, for $g>2$.

\begin{remark}\rm
It is worth mentioning that in \cite{D3}, two alternative bases for $KBSM(H_2)$ are presented in terms of mixed links/ mixed braids (see Figure~\ref{basesall}). 
\begin{figure}[!ht]
\begin{center}
\includegraphics[width=5.2in]{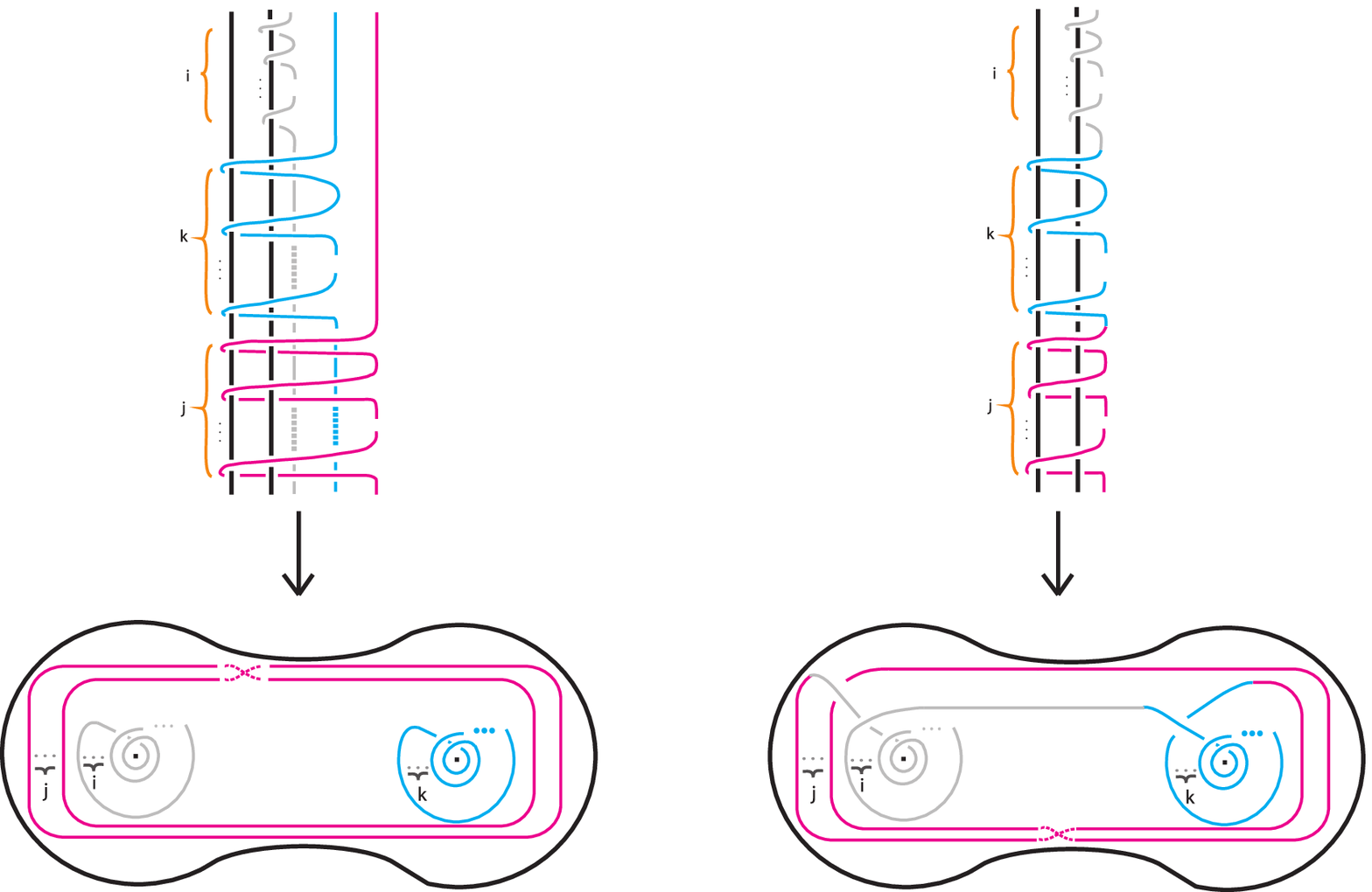}
\end{center}
\caption{Elements in the $B^{\prime}_{H_2}$ and $\mathcal{B}_{H_2}$ bases.}
\label{basesall}
\end{figure}
\end{remark}

\smallbreak

Denote by $L\cup O, L\cup x, L\cup y$ and $L\cup z$ the diagram consisting of the diagram $L$ together with the unknot $O$ and the curves $x, y, z$ respectively. Recall also that for pseudo links, we have to consider oriented diagrams when a pseudo-crossing is resolved using the skein relations. This leads to the following definition:

\begin{defn}\label{pkaufbh2}\rm
Let $L$ be an oriented pseudo link in $H_2$. The {\it pseudo bracket polynomial} of $L$ is defined by means of the relations in Definition~\ref{pkaufb}, together with the following relations:
\begin{figure}[H]
\begin{center}
\includegraphics[width=2.6in]{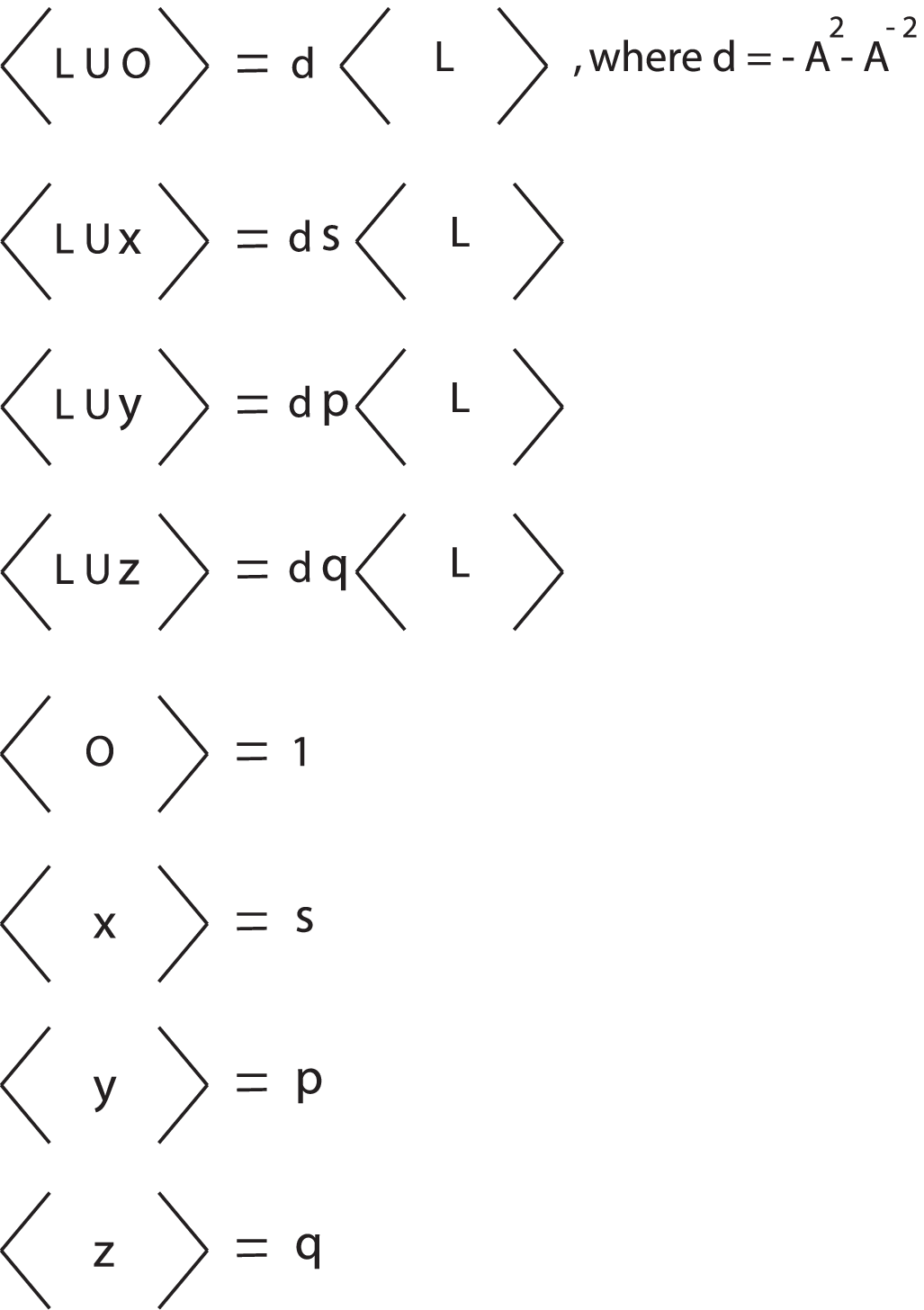}
\end{center}
\label{pkbh22}
\end{figure}
\end{defn}

\begin{figure}[ht]
\begin{center}
\includegraphics[width=4.1in]{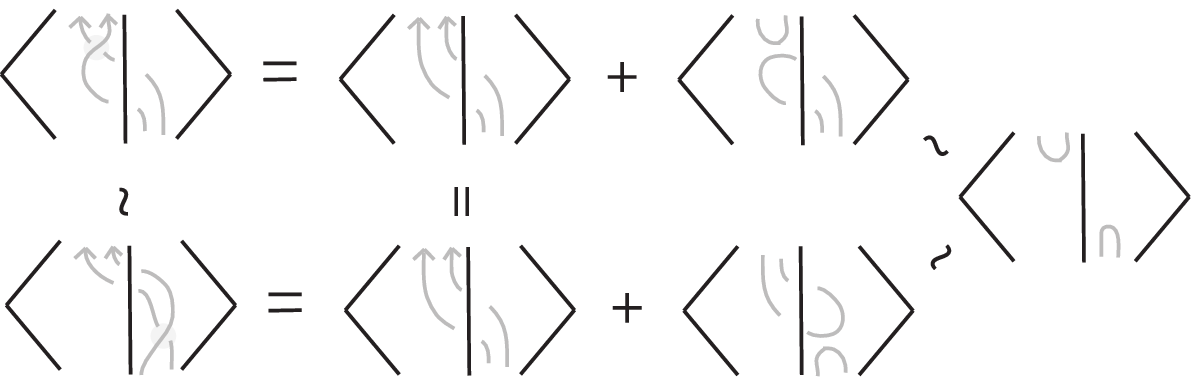}
\end{center}
\caption{The invariance of $<;>$ under the Mixed Reidemeister move 3.}
\label{kf1}
\end{figure}

Obviously, the pseudo bracket polynomial for pseudo links in $H_2$ is invariant under the Reidemeister moves 2 and 3 and under the Mixed Reidemeister move 2. Invariance under the Mixed Reidemeister move 3 is illustrated in Figure~\ref{kf1} (by omitting the coefficients) and invariance under the Mixed pseudo Reidemeister move 3 follows similarly.

Finally, if we consider the normalization of the pseudo bracket polynomial using the writhe, we obtain the {\it normalized pseudo bracket polynomial}, which is an invariant of pseudo links in $H_2$. In particular:

\begin{thm}
Let $L$ be a pseudo diagram of a pseudo link in $H_2$. The polynomial
\[
P_L(A, s, p, q)\ =\ (-A^{-3})^{w(L)}\, <L>,
\]
\noindent where $w(L)$ denotes the writhe of the pseudo link and $<L>$ denotes the pseudo bracket polynomial of $L$, is an invariant of pseudo links in $H_2$. 
\end{thm}

\begin{remark}\rm
A similar result for the case of classical knots in $H_g$ is obtained in \cite{BH}.
\end{remark}

In order to define a polynomial invariant for pseudo links in $H_g$, where $g>2$, we must consider all possible relations of the form $<L \cup w_i>\, =\, (-A^2-A^{-2})\, s_i\, <L>,\, \forall i$, where the $s_i$'s are indeterminate and $w_i$ are basic loops of $H_g$, that is, closed curves that go around the $i^{th}$ pole of $H_g$ once. Indeed, we have the following:

\begin{thm}
Let $L$ be a pseudo diagram of a pseudo link in $H_g$. The polynomial defined by 
\[
P_L(A, w_1,\, ldots,\, w_n)\ =\ (-A^{-3})^{w(L)}\, <L>,
\]
\noindent where $w_i$ denotes all possible non-contractible loops on the $g$-punctured disk, $w(L)$ denotes the writhe of the pseudo link and $<L>$ denotes the pseudo bracket polynomial of $L$, is an invariant of pseudo links in $H_g$. 
\end{thm}

\begin{figure}[ht]
\begin{center}
\includegraphics[width=4.1in]{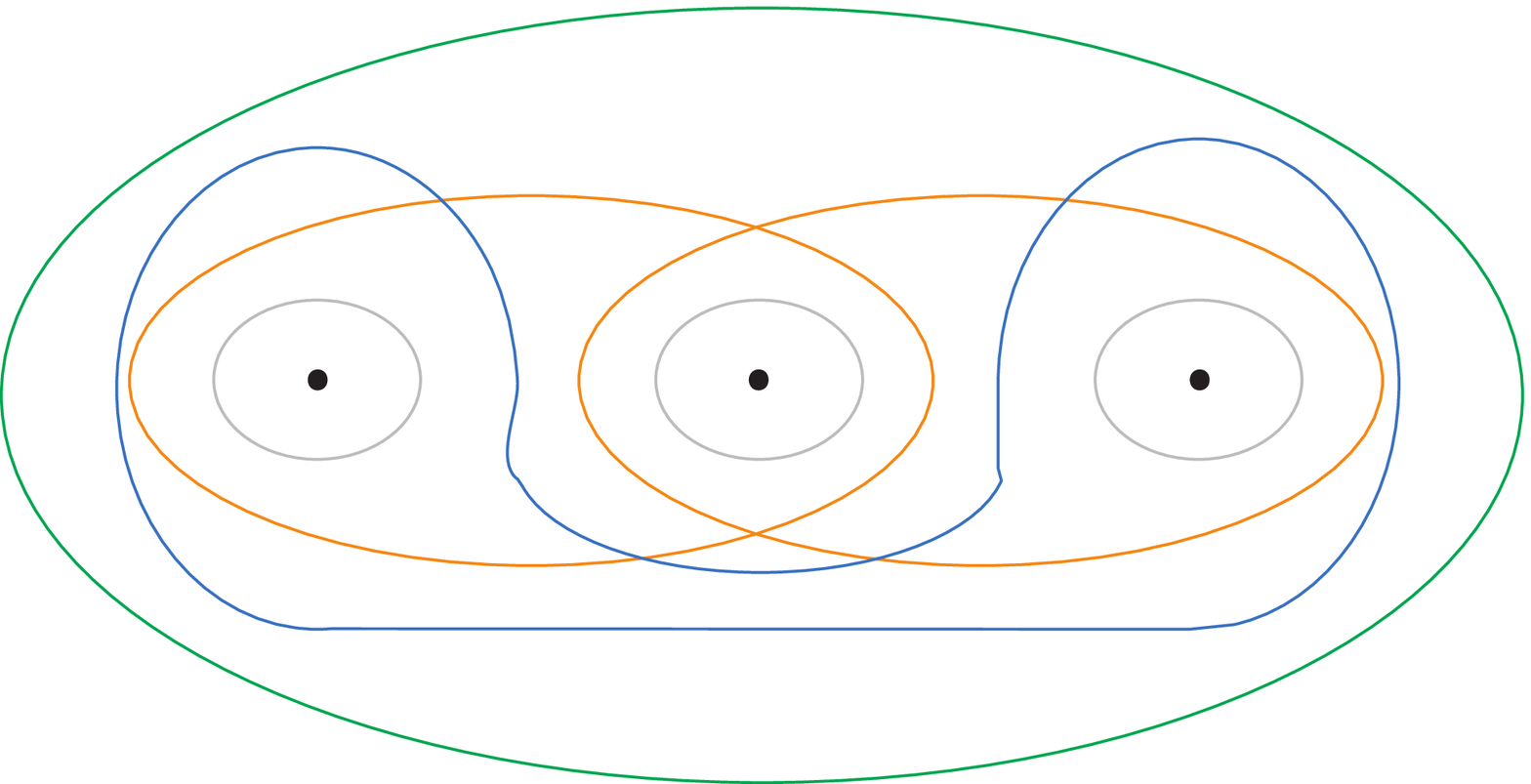}
\end{center}
\caption{Non-contractible loops on the $3$-punctured disk.}
\label{h33}
\end{figure}

For an illustration of non-contractible loops on the $3$-punctured disk, that corresponds to the handlebody of genus $3$, $H_3$, see Figure~\ref{h33}.

\subsection{Pseudo braids in $H_g$}\label{plhg3}

We now introduce the mixed pseudo braids and we present the analogue of the Alexander theorem. For that we need to introduce the notion of {\it geometric mixed pseudo braids} (see also \cite{D5, DL1, La, LR1}).

\begin{figure}[ht]
\begin{center}
\includegraphics[width=5.3in]{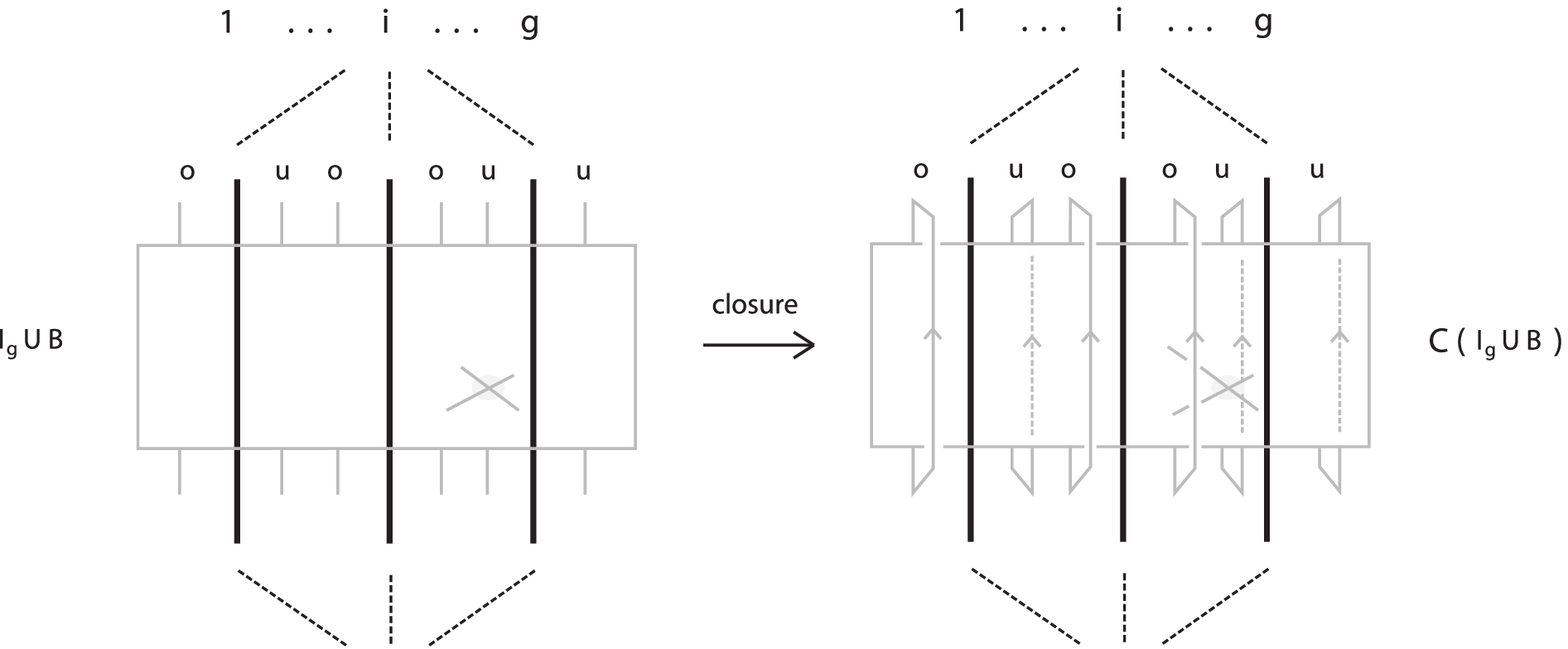}
\end{center}
\caption{An abstract geometric mixed pseudo braid in $H_g$ and its closure.}
\label{bHg}
\end{figure}

\begin{defn}\rm
A \textit{mixed pseudo braid} on $n$ strands, denoted by $I_g\cup B$, is an element of the pseudo braid monoid $PM_{g+n}$ consisting of two disjoint sets of strands, one of which forms the identity braid on $g$-strands, $I_g$, and that represents handlebody of genus $g$, and $n$ strands form the {\it moving subbraid\/} $\beta$ representing the pseudo link $L$ in $H_g$. It is crucial to note that the set of strands that represent $H_g$ have each pair of corresponding endpoints labeled with an $o$ (for over) or $u$ (for under). The reason why we label these strands in that way follows from Definition~\ref{clop}, of the closure operation of geometric mixed braids in $H_g$, and the discussion below. For an illustration see the left hand side of Figure~\ref{bHg}. Moreover, a diagram of a mixed pseudo braid is a braid diagram projected on the plane of $I_g$.
\end{defn}

The main difference between a geometric mixed braid in $H_g$ and a geometric mixed braid in ST lies in the {\it closing} operation that is used to obtain mixed links from geometric mixed braids. In the case of ST, the closing operation is defined in the usual sense (as in the case of classical braids in $S^3$), while in the case of $H_g$ the strands of the fixed part $I_g$ do not participate in the closure operation. In particular, the closure operation in $H_g$ is defined as follows:

\begin{defn}\label{clop}\rm
The {\it closure} $C(I_g \cup B)$ of a geometric mixed pseudo braid $I_g \cup B$ in $H_g$ is realized by joining each pair of corresponding endpoints of the moving part by a vertical segment, either over or under the rest of the braid and according to the label attached to these endpoints (see Figure~\ref{bHg}).
\end{defn}

\begin{remark}\rm
It is crucial to note that different labels on the endpoints of a geometric mixed braid will yield non isotopic links in $H_g$ in general. For more details the reader is referred to \cite{OL}. 
\end{remark}

We now extend the definition of $L$-moves for the case of pseudo links in $H_g$. An $L$-move on a mixed pseudo link is an $L$-move on the moving part of the mixed pseudo link (recall Definition~\ref{lmdefn}). The $L$-moves are similar to the braiding moves involved in the braiding process that we now present.

\smallbreak

\begin{figure}[ht]
\begin{center}
\includegraphics[width=4.3in]{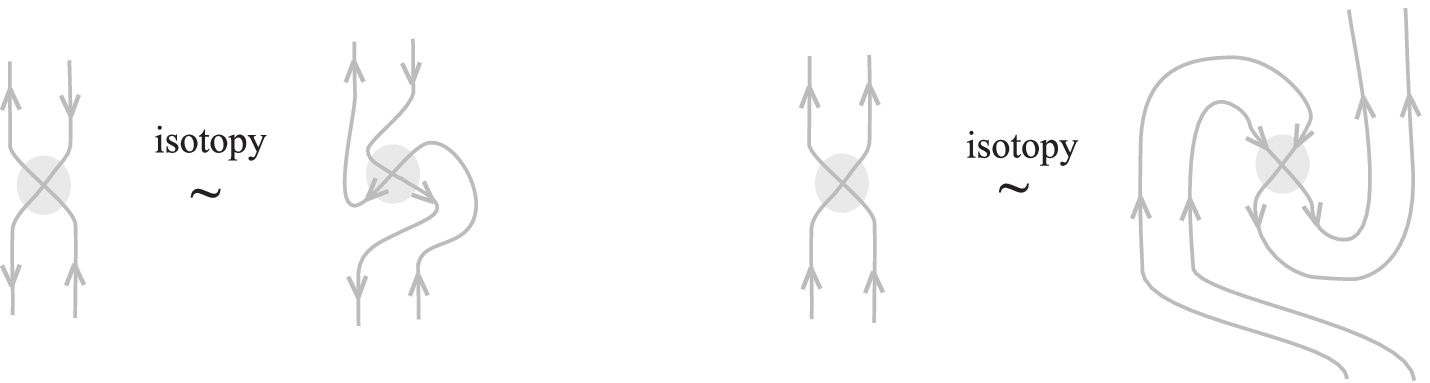}
\end{center}
\caption{Rotating pre-crossings.}
\label{tw}
\end{figure}

In \cite{LR1} a braiding algorithm is presented for classical knots and links in $3$-manifolds (see also \cite{OL} for the case of classical links in $H_g$). The main idea is to keep the arcs of the oriented link diagrams that go downwards with respect to the height function unaffected and replace arcs that go upwards with braid strands. Obviously, these arcs do not belong to the fixed subbraid $I_g$. The same algorithm may be applied for pseudo links in $H_g$, if we first rotate all pseudo-crossings so that all arcs involved will be directed downward as illustrated in Figure~\ref{tw}. Then we may apply the braiding algorithm of \cite{LR1} for the pseudo link (ignoring the pre-crossings).

For the sake of completion, we recall the basic steps of the braiding algorithm of \cite{LR1}:

\bigbreak

\begin{itemize}
\item We first isotope the diagram of the pseudo mixed link as described above. Then, we apply the following braiding algorithm for the pseudo link:
\smallbreak
\item We chose a base-point and we run along the diagram of the mixed pseudo link according to its orientation.
\smallbreak
\item When/If we run along an opposite arc, that is, an arc that goes upwards, we subdivide it into smaller arcs, each containing crossings of one type only as shown in Figure~\ref{upa}. These arcs are called {\it up-arcs}.

\begin{figure}[ht]
\begin{center}
\includegraphics[width=4in]{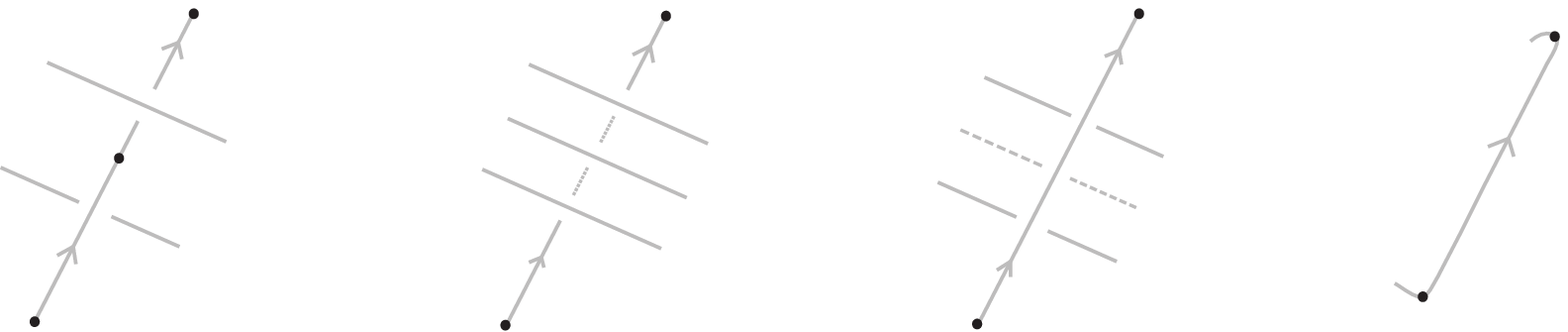}
\end{center}
\caption{Up-arcs.}
\label{upa}
\end{figure}

\item We now label every up-arc with an ``o''or a ``u'', according to the crossings it contains. If it contains no crossings, then the choice is arbitrary.
\smallbreak
\item We perform an $L_o$-move on all up-arcs which were labeled with an ``o'' and an $L_u$-move on all up-arcs which were labeled with an ``u'' (see Figure~\ref{ahg}).

\begin{figure}[ht]
\begin{center}
\includegraphics[width=3.2in]{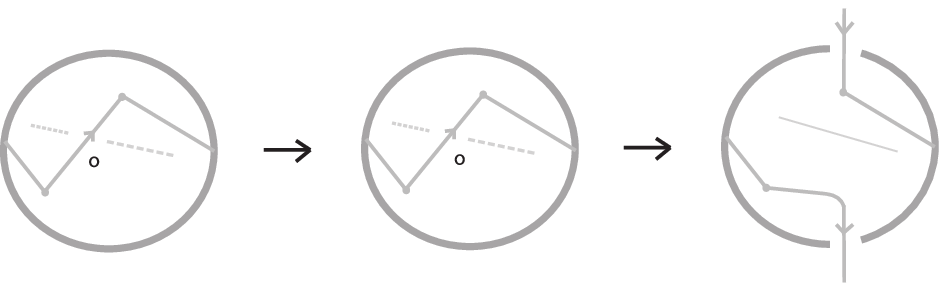}
\end{center}
\caption{Braiding moves for up-arcs.}
\label{ahg}
\end{figure}

\item The result is a geometric mixed pseudo braid whose closure is isotopic to the initial mixed pseudo link.
\end{itemize}

\smallbreak

\noindent Note that the braiding process does not involve the fixed part of the mixed pseudo link.

\smallbreak

The above algorithm provides a proof of the following:

\begin{thm}[{\bf The analogue of the Alexander theorem for pseudo links in $H_g$}]\label{newprpkalex}
Every oriented pseudo link in $H_g$ is isotopic to the closure of a mixed pseudo braid.
\end{thm}

\begin{remark}\rm
Our intention is to obtain invariants for pseudo links in $H_g$ following the Jones approach \cite{Jo} and the results presented in this paper set up the topological background needed. Moreover, we need to obtain the algebraic analogue of the Markov theorem for pseudo links in $H_g$. A geometric analogue is easily obtained by using the $L$-moves (see also \cite{OL} for the case of classical links in $H_g$), but in order to present an algebraic statement for pseudo braid equivalence in $H_g$, we need to introduce and study the {\it mixed pseudo braid monoids}, generated by the standard braiding generators $\sigma_i$, the pseudo-braiding generators $p_i$ and the looping generators $a_j$, illustrated in \ref{lgena}.

\begin{figure}[ht]
\begin{center}
\includegraphics[width=1.7in]{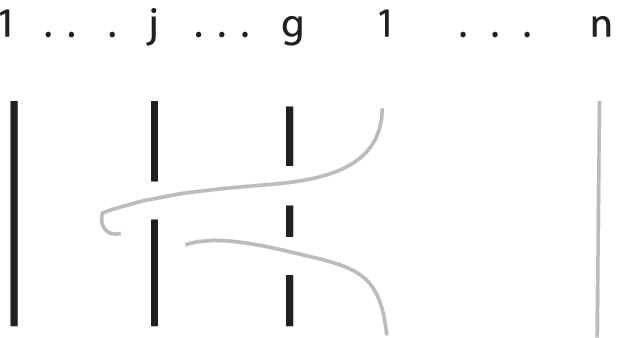}
\end{center}
\caption{ The loop generators $a_j$. }
\label{lgena}
\end{figure}

Obtaining a presentation for mixed pseudo braid monoids is complicated and it will be the subject of a sequel paper. Note that in \cite{V2, D1}, a similar approach for the case of singular links in $H_g$ and for the case of tied links in $H_g$, are presented respectively. Finally, it is worth mentioning that the Jones's approach for obtaining invariants for various knotted objects in various 3-manifolds, has been successfully applied for classical knots in the solid torus (\cite{La1, DL2, D2}), for classical knots in the handlebody of genus 2 (\cite{D3}), for classical knots in the lens spaces $L(p, 1)$ (\cite{DL1, DL3, DL4, DLP, D4}) and for the case of torus knot complements (\cite{D6}).
\end{remark}

\end{document}